\input amstex
\documentstyle{amsppt}

\def\ZZ{\Bbb Z}
\def\QQ{\Bbb Q}
\def\PP{\Bbb P}

\rightheadtext{Arithmetic progressions consisting of unlike powers}
\leftheadtext{N. Bruin, K. Gy\H ory, L. Hajdu and Sz. Tengely}

\topmatter

\title
Arithmetic progressions consisting of unlike powers
\endtitle

\author
N. Bruin$^1$, K. Gy\H ory$^2$, L. Hajdu$^3$ and Sz. Tengely$^4$
\endauthor

\address
\vskip.1cm
\indent
Nils Bruin\hfill\break
\vskip.1cm
\indent
Department of Mathematics\hfill\break
\indent
Simon Fraser University\hfill\break
\indent
Burnaby, BC\hfill\break
\indent
Canada V5A 1S6\hfill\break
\vskip.2cm
\indent
K\'alm\'an Gy\H ory \hfill\break
\indent
Lajos Hajdu \hfill\break
\vskip.1cm
\indent
Number Theory Research Group \hfill\break
\indent
of the Hungarian Academy of Sciences, and \hfill\break
\indent
University of Debrecen \hfill\break
\indent
Institute of Mathematics \hfill\break
\indent
P.O. Box 12 \hfill\break
\indent
4010 Debrecen \hfill\break
\indent
Hungary \hfill\break
\vskip.2cm
\indent
Szabolcs Tengely \hfill\break
\vskip.05cm
\indent
University of Debrecen \hfill\break
\indent
Institute of Mathematics \hfill\break
\indent
P.O. Box 12 \hfill\break
\indent
4010 Debrecen \hfill\break
\indent
Hungary \hfill\break
\vskip.25cm
\endaddress

\email
\hfill\break
\indent
nbruin\@cecm.sfu.ca\hfill\break
\indent
gyory\@math.klte.hu\hfill\break
\indent
hajdul\@math.klte.hu\hfill\break
\indent
tengely\@math.klte.hu
\endemail

\abstract

In this paper we present some new results about unlike powers in
arithmetic progression. We prove among other things that for given $k\geq 4$ and
$L\geq 3$ there are only finitely many arithmetic progressions of the form
$(x_0^{l_0},x_1^{l_1},\ldots,x_{k-1}^{l_{k-1}})$ with $x_i\in{\Bbb Z},$
gcd$(x_0,x_1)=1$ and $2\leq l_i\leq L$ for $i=0,1,\ldots,k-1.$ Furthermore, we
show that, for $L=3$, the progression $(1,1,\ldots,1)$ is the only such progression up to sign.
Our proofs involve some well-known theorems of Faltings \cite{F}, Darmon and
Granville \cite{DG} as well as Chabauty's method applied to superelliptic
curves.

\endabstract

\endtopmatter

\document

\footnotemark"{ }"\footnotetext"{ }"{\hskip-.1cm {\it 2000 Mathematics
Subject Classification:} 11D41.\hfill\break
\noindent
$^1$Research supported in part by National Science and Engineering Research
Council Canada (NSERC).\hfill\break
\noindent
$^2$Research supported in part by grants T42985 and T38225 of the Hungarian National
Foundation for Scientific Research (HNFSR).\hfill\break
\noindent
$^3$Research supported in part by grants T42985 and T48791 of the HNFSR
and by the J\'anos Bolyai Research Fellowship of the Hungarian Academy of Sciences.\hfill\break
\noindent
$^4$Research supported in part by grant T48791 of the HNFSR.
}

\heading{1. Introduction}
\endheading

By a classical result of Euler, which apparently was already known to Fermat
(see \cite{D} pp. 440 and 635), four distinct squares cannot form an arithmetic
progression. Darmon and Merel \cite{DM} proved that, apart from trivial cases,
there do not exist 3-term arithmetic progressions consisting of $l$-th powers,
provided $l\geq 3.$ More generally, perfect powers from products of consecutive
terms in arithmetic progression have been extensively studied in a great number
of papers; see e.g. \cite{T}, \cite{Sh} and \cite{BBGyH} and the references
there. In our article we deal with the following problem.

\proclaim{Question} For all $k\geq 3$
characterize the non-constant arithmetic
progressions $$(h_0,h_1,\ldots,h_{k-1})$$ with gcd$(h_0,h_1)=1$ such that each
$h_i=x_i^{l_i}$ for some $x_i\in{\Bbb Z}$ and $l_i\geq 2$.
\endproclaim

Note that we impose the seemingly artificial primitivity condition
gcd$(h_0,h_1)=1$. In
case the $h_i$ are all like powers, the homogeneity of the conditions ensures
that up to scaling, we can assume gcd$(h_0,h_1)=1$ without loss of generality.
If we do not take all $l_i$ equal, however, there are infinite families that are
not quite trivial, but are characterized by the fact they have a fairly large
common factor in their terms; see the examples below Theorem 3.

By a recent result of Hajdu \cite{H} the ABC conjecture implies that if
$$(x_0^{l_0},x_1^{l_1},\ldots,x_{k-1}^{l_{k-1}})$$ is an arithmetic progression
with gcd$(x_0,x_1)=1$ and $l_i\geq 2$ for each $i,$ then $k$ and the $l_i$ are
bounded.
Furthermore, he shows unconditionally that $k$ can be bounded above in terms of
$\max_i\{l_i\}.$ In fact Hajdu proves these results for more general arithmetic
progressions which satisfy the assumptions (i), (ii) of our Theorem 2 below.

As is known (see e.g. \cite{M},\cite{DG},\cite{PT},\cite{T1},\cite{T2} and the
references given there), there exist integers $l_0,l_1,l_2\geq 2$ for which there are infinitely many primitive arithmetic progressions of the form $(x_0^{l_0},x_1^{l_1},x_2^{l_2}).$  In these progressions the exponents in question always satisfy the condition
$$
\frac{1}{l_0}+\frac{1}{l_1}+\frac{1}{l_2}\geq 1.
$$
One would, however, expect only very few primitive arithmetic progressions of length at
least four and consisting entirely from powers at least two. A definitive answer
to the above question seems beyond present techniques. As in \cite{H}, we restrict
the size of the exponents $l_i$ and prove the following finiteness result:

\proclaim{Theorem 1} Let $k\geq 4$ and $L\geq 2$. There are only finitely many $k$-term integral
arithmetic progressions $(h_0,h_1,\ldots,h_{k-1})$ such that $\gcd(h_0,h_1)=1$ and  $h_i=x_i^{l_i}$  with some $x_i\in{\Bbb Z}$ and $2\leq l_i\leq L$ for $i=0,1,\ldots,k-1.$
\endproclaim

The proof of this theorem uses that for each of the finitely many possible
exponent vectors $(l_0,\ldots,l_{k-1})$, the primitive arithmetic progressions of the
form $(x_0^{l_0},\ldots,x_{k-1}^{l_{k-1}})$ correspond to the rational points on finitely
many algebraic curves. In most cases, these curves are of genus larger than $1$
and thus, by Faltings' theorem \cite{F}, give rise to only finitely many solutions.

In fact, our Theorem 1 above is a direct consequence of the following more general
result and a theorem by Euler on squares in arithmetic progression.
For a finite set of primes
$S$, we write ${\Bbb Z}_S^*$ for the set of rational integers not divisible by primes outside $S.$

\proclaim{Theorem 2}
Let $L,k$ and $D$ be positive integers with $L\geq 2,k\geq 3,$ and let $S$ be a finite set of primes. Then there are at most
finitely many arithmetic progressions $(h_0,h_1,\ldots,h_{k-1})$ 
satisfying the following conditions:
\roster
\item"(i)" For $i=0,\ldots,k-1$, there exist $x_i\in\Bbb{Z}$, $2\leq l_i\leq L$ and
$\eta_i\in\Bbb{Z}_S^*$ such that $$h_i=\eta_i\,x_i^{l_i},$$
\item"(ii)" $\gcd(h_0,h_1)\leq D,$
\item"(iii)" either $k\geq 5$, or $k=4$ and $l_i\geq 3$ for some $i,$ or $k=3$ and $\frac{1}{l_0}+\frac{1}{l_1}+\frac{1}{l_2}<1.$
\endroster
\endproclaim

\noindent
{\bf Remark.}
In (iii) the assumptions concerning the exponents $l_i$ are necessary. For $k=3$ this was seen above. In case of $k=4$ the condition $l_i\geq 3$ for some $i$ cannot be omitted as is shown by e.g. the arithmetic progression $x_0^2,x_1^2,x_2^2,73x_3^2$ with $S=\{73\}.$ We have the homogeneous system of equations
$$
\align
x_0^2+x_2^2&=2x_1^2\\
x_1^2+73x_3^2&=2x_2^2.
\endalign
$$
A non-singular intersection of two quadrics in $\Bbb{P}^3$ is a genus 1 curve. If
there is a rational point on it, it is isomorphic to its Jacobian - an
elliptic curve. In this example the elliptic curve has infinitely many
rational points. Therefore we also have infinitely many rational solutions
$(x_0:x_1:x_2:x_3).$ After rescaling, those all give primitive integral solutions as well.

For small $l_i$ we can explicitly find the parametrizing algebraic
curves and, using Chabauty's method, the rational points on them. This allows us to prove:

\proclaim{Theorem 3}
Let $k\geq 4,$ and suppose that $(h_0,h_1,\ldots,h_{k-1})=(x_0^{l_0},x_1^{l_1},\ldots,x_{k-1}^{l_{k-1}})$ is a
primitive integral arithmetic progression with $x_i\in{\Bbb Z}$ and $2\leq l_i\leq
3$ for $i=0,1,\ldots,k-1$. Then $$(h_0,h_1,\ldots,h_{k-1})=\pm(1,1,\ldots,1).$$
\endproclaim

The proof is rather computational in nature and uses $p$-adic methods to derive
sharp bounds on the number of rational points on specific curves. The methods
are by now well-established. Of particular interest to the connoisseur would be
the argument for the curve ${\Cal C}_4$ in Section 3, where we derive that an
elliptic curve has rank $0$ and a non-trivial Tate-Shafarevich group by doing a
full $2$-descent on an isogenous curve and the determination of the solutions to
equation \thetag{7}. The novelty for the latter case lies in the fact that, rather than
considering a hyperelliptic curve, we consider a superelliptic curve of the form
$$f(x)= y^3,\text{ with }\deg(f)=6.$$
We then proceed similarly to \cite{B}. We determine an extension $K$ over which
$f(x)=g(x)\cdot h(x)$, with $g,h$ both cubic. We then determine that $\QQ$-rational
solutions to $f(x)=y^3$ by determining, for finitely many values $\delta$, the
$K$-rational points on the genus $1$ curve $g(x)=\delta y_1^3$, with $x\in\QQ$.

\medskip
\noindent
{\bf Remark.} The condition gcd$(h_0,h_1)=1$ in Theorems 1 and 3 is necessary. This
can be illustrated by the following examples with $k=4.$ Note
that the progressions below can be ``reversed'' to get examples for the
opposite orders of the exponents $l_0,l_1,l_2,l_3$.

$\bullet$ In case of $(l_0,l_1,l_2,l_3)=(2,2,2,3)$
$$
((u^2-2uv-v^2)f(u,v))^2,((u^2+v^2)f(u,v))^2,((u^2+2uv-v^2)f(u,v))^2,(f(u,v))^3
$$
is an arithmetic progression for any $u,v\in{\Bbb Z}$, where
$f(u,v)=u^4+8u^3v+2u^2v^2-8uv^3+v^4$.

$\bullet$ In case of $(l_0,l_1,l_2,l_3)=(2,2,3,2)$
$$
((u^2-2uv-2v^2)g(u,v))^2,((u^2+2v^2)g(u,v))^2,(g(u,v))^3,((u^2+4uv-2v^2)g(u,v))^2
$$
is an arithmetic progression for any $u,v\in{\Bbb Z}$, where
$g(u,v)=u^4+4u^3v+8u^2v^2-8uv^3+4v^4$.

\heading{2. Auxiliary results}
\endheading

The proof of Theorem 2 depends on the following well-known result by Darmon
and Granville \cite{DG}.

\proclaim{Theorem A} 
Let $A,B,C$ and $r,s,t$ be non-zero integers with $r,s,t\geq 2,$ and let $S$ be a finite set of primes.
Then there exists a number field $K$ such that all solutions $x,y,z\in{\Bbb Z}$
with gcd$(x,y,z)\in{\Bbb Z}_S^*$ to the equation
$$Ax^r+By^s=Cz^t$$
correspond, up to weighted projective equivalence, to $K$-rational points on
some algebraic curve $X_{r,s,t}$ defined over $K$. Putting
$u=-Ax^r/Cz^t$, the curve $X$ is a Galois-cover of the $u$-line of degree $d$,
unramified outside $u\in\{0,1,\infty\}$ and with ramification indices
$e_0=r,e_1=s,e_2=t$.
Writing $\chi(r,s,t)=1/r+1/s+1/t$ and $g$ for the genus of $X$, we find
\roster
\item"$\bullet$" if $\chi(r,s,t)>1$ then $g=0$ and $d=2/\chi(r,s,t)$,
\item"$\bullet$" if $\chi(r,s,t)=1$ then $g=1$,
\item"$\bullet$" if $\chi(r,s,t)<1$ then $g>1$.
\endroster
\endproclaim

The two results below will be useful for handling special progressions,
containing powers with small exponents. The first one deals with the
quadratic case.

\proclaim{Theorem B} Four distinct squares cannot form an arithmetic
progression.
\endproclaim

\demo{Proof} The statement is a simple consequence of a classical result of
Euler (cf. \cite{M}, p. 21), which was already known by Fermat (see \cite{D}
pp. 440 and 635).
\qed
\enddemo

We also need a classical result on a cubic equation.

\proclaim{Theorem C} The equation $x^3+y^3=2z^3$ has the only solutions
$(x,y,z)=\pm(1,1,1)$ in non-zero integers $x,y,z$ with gcd$(x,y,z)=1$.
\endproclaim

\demo{Proof} See Theorem 3 in \cite{M} on p. 126.
\qed
\enddemo

The next lemma provides the parametrization of the solutions of certain
ternary Diophantine equations.

\proclaim{Lemma} All solutions of the equations
$$
\text{\rom{i)}}\ 2b^2-a^2=c^3, \ \ \ \text{\rom{ii)}}\ a^2+b^2=2c^3,\ \ \ \text{\rom{iii)}}\ a^2+2b^2=3c^3,
\ \ \ \text{\rom{iv)}}\ 3b^2-a^2=2c^3,
$$
$$
\text{\rom{v)}}\ 3b^2-2a^2=c^3, \ \ \ \text{\rom{vi)}}\ a^2+b^2=2c^2,\ \ \ \text{\rom{vii)}}\ 2a^2+b^2=3c^2,
\ \ \ \text{\rom{viii)}}\ a^2+3b^2=c^2
$$
in integers $a$, $b$ and $c$ with gcd$(a,b,c)=1$ are given by the following
parametrizations:

$$\halign to \hsize{#\tabskip=8pt plus.15em\hfil&
   #\tabskip=8pt plus.35em\hfil&
   #\tabskip=8pt plus.35em\hfil&
#\tabskip=0pt\hfil\cr
\text{\rom{i)}}&$a=\pm(x^3+6xy^2)$&or&$a=\pm(x^3+6x^2y+6xy^2+4y^3)$\cr
&$b=\pm(3x^2y+2y^3)$&&$b=\pm(x^3+3x^2y+6xy^2+2y^3)$\cr

\text{\rom{ii)}}&$a=\pm(x^3-3x^2y-3xy^2+y^3)$&&\cr
&$b=\pm(x^3+3x^2y-3xy^2-y^3)$&&\cr
\text{\rom{iii)}}&$a=\pm(x^3-6x^2y-6xy^2+4y^3)$&&\cr
&$b=\pm(x^3+3x^2y-6xy^2-2y^3)$&&\cr
\text{\rom{iv)}}&$a=\pm(x^3+9x^2y+9xy^2+9y^3)$&or&$a=\pm(5x^3+27x^2y+45xy^2+27y^3)$\cr
&$b=\pm(x^3+3x^2y+9xy^2+3y^3)$&&$b=\pm(3x^3+15x^2y+27xy^2+15y^3)$\cr
\text{\rom{v)}}&$a=\pm(x^3+9x^2y+18xy^2+18y^3)$&or&$a=\pm(11x^3+81x^2y+198xy^2+162y^3)$\cr
&$b=\pm(x^3+6x^2y+18xy^2+12y^3)$&&$b=\pm(9x^3+66x^2y+162xy^2+132y^3)$\cr
\text{\rom{vi)}}&$a=\pm(x^2-2xy-y^2)$&&\cr
&$b=\pm(x^2+2xy-y^2)$&&\cr
\text{\rom{vii)}}&$a=\pm(x^2+2xy-2y^2)$&&\cr
&$b=\pm(x^2-4xy-2y^2)$&&\cr
\text{\rom{viii)}}&$a=\pm(x^2-3y^2)/2$&&\cr
&$b=\pm xy$&&\cr
}
$$
\noindent
Here $x$ and $y$ are coprime integers and the $\pm$ signs can
be chosen independently.
\endproclaim

\demo{Proof} The statement can be proved via factorizing the
expressions in the appropriate number fields. More precisely, we have to
work in the rings of integers of the following fields:
$\Bbb{Q}(\sqrt{-2}),\Bbb{Q}(i),\Bbb{Q}(\sqrt{2}),\Bbb{Q}(\sqrt{3}),
\Bbb{Q}(\sqrt{6})$. Note that the class number is
one in all of these fields. As the method of the proof of the separate cases
are rather similar, we give it only in two characteristic instances,
namely for the cases i) and vii).

\vskip.3cm

\noindent
i) In ${\Bbb Z}[\sqrt{2}]$ we have
$$
(a+\sqrt{2}b)(a-\sqrt{2}b)={(-c)}^3.
$$
Using gcd$(a,b)=1$, a simple calculation gives that
$$
\text{gcd}(a+\sqrt{2}b,a-\sqrt{2}b)\mid 2\sqrt{2}
$$
in ${\Bbb Z}[\sqrt{2}]$. Moreover, $1+\sqrt{2}$ is a fundamental unit of
${\Bbb Z}[\sqrt{2}]$, and the only roots of unity are $\pm 1$, which are
perfect cubes. Hence we have
$$
a+\sqrt{2}b={(1+\sqrt{2})}^\alpha {(\sqrt{2})}^\beta
{(x+\sqrt{2}y)}^3,\tag{1}
$$
where $\alpha\in\{-1,0,1\}$, $\beta\in\{0,1,2\}$ and $x,y$ are some rational
integers. By taking norms, we immediately obtain that $\beta=0$. If
$\alpha=0$, then expanding the right hand side of \thetag{1} we get
$$
a=x^3+6xy^2,\ \ \ b=3x^2y+2y^3.
$$
Otherwise, when $\alpha=\pm 1$ then \thetag{1} yields
$$
a=x^3\pm 6x^2y+6xy^2\pm 4y^3,\ \ \ b=\pm x^3+3x^2y\pm 6xy^2+2y^3.
$$
In both cases, substituting $-x$ and $-y$ for $x$ and $y$,
respectively, we obtain the parametrizations given in the statement.
Furthermore, observe that the coprimality of $a$ and $b$ implies gcd$(x,y)=1$.

\vskip.25cm

\noindent
vii) By factorizing in ${\Bbb Z}[\sqrt{-2}]$ we obtain
$$
(b+\sqrt{-2}a)(b-\sqrt{-2}a)=3c^2.
$$
Again, gcd$(a,b)=1$ implies that
$$
\text{gcd}(b+\sqrt{-2}a,b-\sqrt{-2}a)\mid 2\sqrt{-2}
$$
in ${\Bbb Z}[\sqrt{-2}]$. Note that ${\Bbb Z}[\sqrt{-2}]$ has no other units than
$\pm 1.$ Since $2=-(\sqrt{-2})^2,$ we can write
$$
b+\sqrt{-2}a={(-1)}^\alpha {(1+\sqrt{-2})}^\beta {(1-\sqrt{-2})}^\gamma
{(\sqrt{-2})}^\delta {(x+\sqrt{-2}y)}^2,\tag{2}
$$
where $\alpha,\beta,\gamma,\delta\in\{0,1\}$ and $x,y$ are some rational integers.
By taking norms, we immediately get that $\delta=0$ and $\beta+\gamma=1$.
In these cases, by expanding the right hand side of \thetag{2} we obtain
(choosing the $\pm$ signs appropriately) that
$$
a=\pm(\pm x^2+2xy\mp y^2),\ \ \ b=\pm(x^2\mp 4xy-2y^2).
$$
Substituting $-x$ and $-y$ in places of $x$ and $y$, respectively, we get
the parametrizations indicated in the statement. Again, gcd$(a,b)=1$ gives
gcd$(x,y)=1$.
\qed
\enddemo

\heading{3. Proofs of the Theorems}
\endheading
Note that Theorem 1 directly follows from Theorem B and Theorem 2. Hence we begin with the proof of the latter statement.
\demo{Proof of Theorem 2}
Since an arithmetic progression of length $k>5$ contains an arithmetic
progression of length $5$, we only have to consider the cases $k=5,4$ and $3$. The
condition that $2\leq l_i\leq L$ leaves only finitely many possibilities for the
exponent vector $\underline{l}=(l_0,\ldots,l_{k-1})$. Therefore, it suffices to prove the
finiteness for a given exponent vector $\underline{l}$.

Note that if $h_i=\eta_i x_i^{l_i}$ for some $\eta_i\in{\Bbb Z}_S^*$, then without
loss of generality, $\eta_i$ can be taken to be $l_i$-th power free. This means
that, given $\underline{l}$, we only need to consider finitely many vectors
$\underline{\eta}=(\eta_0,\ldots,\eta_{k-1})$. Hence, we only need to prove the theorem
for $k=3,4,5$, and $\underline{l}$ and $\underline{\eta}$ fixed. Note that if
$\gcd(h_0,h_1)\leq D$,
then certainly $\gcd(x_i,x_j) \leq D$. 
We enlarge $S$ with all primes up to $D$.

We write $n=h_1-h_0$ for the increment of the arithmetic progression. With
$k,\underline{l},\underline{\eta}$ fixed, the theorem will be proved if we show that the following
system of equations has only finitely many solutions:
\roster
\item"(a)" $\eta_i x_i^{l_i}-\eta_j x_j^{l_j}=(i-j)n$ for all $0\leq i<j\leq k-1$.
\item"(b)" $(x_0,\ldots,x_{k-1})\in\Bbb{Z}^k$ with gcd$(x_0,x_1)\leq D.$
\endroster
Hence, we need to solve
$$(j-m)\eta_ix_i^{l_i}+(m-i)\eta_jx_j^{l_j}+(i-j)\eta_mx_m^{l_m}=0\text{ for all }
0\leq m,i,j\leq k-1.$$
For $m=0,i=1$, we obtain that each of our solutions would give rise to a
solution to
$$
j\eta_1x_1^{l_1}-\eta_jx_j^{l_j}+(1-j)\eta_0x_0^{l_0}=0. \tag{3}
$$
By applying Theorem A we see that such solutions give rise to
$K_j$-rational points on some algebraic curve $C_j$ over some number field $K_j$.
Furthermore, putting
$$u=\frac{\eta_1x_1^{l_1}}{\eta_0x_0^{l_0}},$$
we obtain that
$C_j$ is a Galois-cover of the $u$-line, with ramification indices $l_0,l_1,l_j$
over $u=\infty,0,j/(j-1)$ respectively and unramified elsewhere.

If $k=3$, we recover the approach of Darmon and Granville. Theorem A immediately
implies that if $1/l_0+1/l_1+1/l_2<1$ then $C_2$ has genus larger than $1$ and
thus (by Faltings) has only finitely many rational points. This establishes the
desired finiteness result.

If $k=4$, we are interested in solutions to \thetag{3}
for $j=2,3$
simultaneously. Let $M$ be a number field containing both
$K_2$ and $K_3$. Then the solutions we are interested in,
correspond to $M$-rational points on $C_2$ and $C_3$ that give rise to the same
value of $u$, i.e., we want the rational points on the fibre product
$C_2\times_{u} C_3$. This fibre product is again Galois and has ramification
indices at least $l_0,l_1,l_2,l_3$ over $u=\infty,0,2,\frac{3}{2},$ respectively. Since $C_2\times_{u}
C_3$ is Galois over the $u$-line, all its connected components have the same
genus and degree, say, $d$. Writing $g$ for the genus of this component,
the Riemann-Hurwitz formula gives us
$$2g-2\geq
d\left(2-\frac{1}{l_0}-\frac{1}{l_1}-\frac{1}{l_2}-\frac{1}{l_3}\right).$$
Hence, we see that $g\leq 1$ only if $l_0=l_1=l_2=l_3=2.$ For
other situations, we have $g\geq 2$, so $C_2\times_u C_3$ has only finitely many
$M$-rational points.

If $k=5$, we argue similarly, but now we consider $C_2\times_u
C_3\times_u C_4$, with ramification indices at least $l_0,l_1,l_2,l_3,l_4$
over $u=0,\infty,1,\frac{3}{2},\frac{4}{3},$ respectively. Hence, we obtain
$$2g-2\geq
d\left(3-\frac{1}{l_0}-\frac{1}{l_1}-\frac{1}{l_2}-\frac{1}{l_3}-\frac{1}{l_4}\right),$$
so that $g \geq 2$ in all cases. 

\qed
\enddemo

\demo{Proof of Theorem 3} The proof involves some explicit computations
that are too involved to do either by hand or reproduce here on paper. Since
the computations are by now completely standard, we choose not to bore the
reader with excessive details and only give a conceptual outline of the
proof. For full details, we refer the reader to the electronic resource
\cite{notes}, where a full transcript of a session using the computer
algebra system MAGMA \cite{magma} can be found. We are greatly indebted
to all contributors to this system. Without their work, the computations
sketched here would not at all have been trivial to complete.

It suffices to prove the assertion for $k=4.$ We divide the proof into several parts, according to the exponents of the
powers in the arithmetic progression. If $(l_0,l_1,l_2,l_3)=(2,2,2,2)$,
$(3,3,3,3)$, $(2,3,3,3)$ or $(3,3,3,2)$, then our
statement follows from Theorems B and C. We handle the remaining cases by
Chabauty's method. We start with those cases where the classical variant
works. After that we consider the cases where we have to resort to considering
some covers of elliptic curves.

\vskip.15cm

\noindent
{\bf The cases $(l_0,l_1,l_2,l_3)=(2,2,2,3)$ and $(3,2,2,2)$.}

\vskip.15cm

\null From the method of our proof it will be clear that by symmetry we may
suppose $(l_0,l_1,l_2,l_3)=(2,2,2,3)$. That is, the progression is of
the form $x_0^2,x_1^2,x_2^2,x_3^3$. Applying part i) of our Lemma to the last
three terms of the progression, we get that either
$$
x_1=\pm(x^3+6xy^2),\ \ \ x_2=\pm(3x^2y+2y^3)
$$
or
$$
x_1=\pm(x^3+6x^2y+6xy^2+4y^3),\ \ \ x_2=\pm(x^3+3x^2y+6xy^2+2y^3)
$$
where $x,y$ are some coprime integers in both cases.

In the first case by $x_0^2=2x_1^2-x_2^2$ we get
$$
x_0^2=2x^6+15x^4y^2+60x^2y^4-4y^6.
$$
Observe that $x\neq 0$. By putting $Y=x_0/x^3$ and $X=y^2/x^2$ we obtain
the elliptic equation
$$
Y^2=-4X^3+60X^2+15X+2.
$$ 
A straightforward calculation with MAGMA gives that the elliptic curve described
by this equation has no affine rational points.

In the second case by the same assertion we obtain
$$
x_0^2=x^6+18x^5y+75x^4y^2+120x^3y^3+120x^2y^4+72xy^5+28y^6.
$$
If $y=0$, then the coprimality of $x$ and $y$ yields $x=\pm 1$, and we
get the trivial progression $1,1,1,1$. So assume that $y\neq 0$ and let
$Y=x_0/y^3$, $X=x/y$. By these substitutions we are led to the
hyperelliptic (genus two) equation
$$
{\Cal C}_1:Y^2=X^6+18X^5+75X^4+120X^3+120X^2+72X+28.
$$
We show that ${\Cal C}_1({\Bbb Q})$ consists only of the two points on 
${\Cal C}_1$ above $X=\infty,$ denoted by $\infty^+$ and $\infty^-$.

The order of ${\Cal J}_{\text{tors}}({\Bbb Q})$ (the torsion subgroup of
the Mordell-Weil group ${\Cal J}({\Bbb Q})$ of the Jacobian of ${\Cal C}_1$) is a divisor of
$\gcd(\#{\Cal J}({\Bbb F}_5),\#{\Cal J}({\Bbb F}_7))=\gcd(21,52)=1$. Therefore
the torsion subgroup is trivial. Moreover, using the algorithm of M. Stoll
\cite{St} implemented in MAGMA we get that the rank of
${\Cal J}({\Bbb Q})$ is at most one. As the divisor $D=[\infty^+-\infty^-]$
has infinite order, the rank is exactly one. Since the rank
of ${\Cal J}({\Bbb Q})$ is less than the genus of ${\Cal C}_1$, we can
apply Chabauty's method \cite{C} to obtain a bound for the number of
rational points on ${\Cal C}_1$. For applications of the method on
related problems, we refer to \cite{CF}, \cite{Fl}, \cite{FPS}, \cite{P}.

As the rank of ${\Cal J}({\Bbb Q})$ is one and the torsion is trivial, we
have ${\Cal J}({\Bbb Q})=\langle D_0 \rangle$ for some
$D_0\in{\Cal J}({\Bbb Q})$ of infinite order. A simple computation
$\pmod{13}$ shows that $D\notin 5{\Cal J}({\Bbb Q})$, and a similar
computation $\pmod{139}$ yields that $D\notin 29{\Cal J}({\Bbb Q})$. Hence
$D=kD_0$ with $5\nmid k$, $29\nmid k$. The reduction of ${\Cal C}_1$ modulo
$p$ is a curve of genus two for any prime $p\neq 2,3.$ We take
$p=29$. Using Chabauty's method as implemented in MAGMA by Stoll, we find
that there are at most two rational points on ${\Cal C}_1$. Therefore we
conclude that ${\Cal C}_1(\Bbb{Q})=\{\infty^+,\infty^-\}$, which proves
the theorem in this case.

\vskip.15cm

\noindent
{\bf The cases $(l_0,l_1,l_2,l_3)=(2,2,3,2)$ and $(2,3,2,2)$.} 

\vskip.15cm

Again, by symmetry we may suppose that $(l_0,l_1,l_2,l_3)=(2,2,3,2)$.
Then the progression is given by $x_0^2,x_1^2,x_2^3,x_3^2$. Now from part
iii) of our Lemma, applied to the terms with indices $0,2,3$ of the
progression, we get
$$
x_0=\pm(x^3-6x^2y-6xy^2+4y^3),\ \ \ x_3=\pm(x^3+3x^2y-6xy^2-2y^3)
$$
where $x,y$ are some coprime integers. Using $x_1^2=(2x_0^2+x_3^2)/3$
we obtain
$$
x_1^2=x^6-6x^5y+15x^4y^2+40x^3y^3-24xy^5+12y^6.
$$
If $y=0$, then in the same way as before we deduce that the only
possibility is given by the progression $1,1,1,1$. Otherwise, if
$y\neq 0$, then write $Y=x_1/y^3$, $X=x/y$ to get the hyperelliptic
(genus two) curve
$$
{\Cal C}_2:Y^2=X^6-6X^5+15X^4+40X^3-24X+12.
$$
By a calculation similar to that applied in the previous case (but now
with $p=11$ in place of $p=29$) we get that ${\Cal C}_2(\Bbb{Q})$
consists only of the points $\infty^+$ and $\infty^-$. Hence the
statement is proved also in this case.

\vskip.15cm

\noindent {\bf The cases $(l_0,l_1,l_2,l_3)=(3,2,3,2)$ and $(2,3,2,3)$.}

\vskip.15cm

As before, without loss of generality we may assume
$(l_0,l_1,l_2,l_3)=(3,2,3,2)$. Then the progression is given by
$x_0^3,x_1^2,x_2^3,x_3^2$. We have
$$
x_1^2=\frac{x_0^3+x_2^3}{2},\ \ \ x_3^2=\frac{-x_0^3+3x_2^3}{2}.\tag{4}
$$
We note that $x_2=0$ implies $x_1^2=-x_3^2$, hence $x_1=x_3=0$. So we may
assume that $x_2\neq 0$, whence we obtain from \thetag{4} that
$$
\left(\frac{2x_1x_3}{x_2^3}\right)^2=-\left(\frac{x_0}{x_2}\right)^6+
2\left(\frac{x_0}{x_2}\right)^3+3.
$$
Thus putting $Y=2x_1x_3/x_2^3$ and $X=x_0/x_2$, it is sufficient to
find all rational points on the hyperelliptic curve
$$
{\Cal C}_3:Y^2=-X^6+2X^3+3.
$$
We show that ${\Cal C}_3({\Bbb Q})=\{(-1,0),(1,\pm 2)\}$.

Using MAGMA we obtain that the rank of the Jacobian ${\Cal J}({\Bbb Q})$ of
${\Cal C}_3(\Bbb Q)$ is at most one, and the torsion subgroup
${\Cal J}_{\text{tors}}({\Bbb Q})$ of ${\Cal J}({\Bbb Q})$ consists of the
elements ${\Cal O}$ and
$[(\frac{1-\sqrt{3}i}{2},0)+(\frac{1+\sqrt{3}i}{2},0)-\infty^+-\infty^-]$.
As the divisor $D=[(-1,0)+(1,-2)-\infty^+-\infty^-]$ has infinite order, the
rank of ${\Cal J}({\Bbb Q})$ is exactly one. The only Weierstrass point on
${\Cal C}_3$ is $(-1,0)$. We proceed as before, using the primes $7$ and $11$ 
in this case. We
conclude that $(1,\pm 2)$ are the only non-Weierstrass points on
${\Cal C}_3$. It is easy to check that these points give rise only to the
trivial arithmetic progression, so our theorem is proved also in this case.

\vskip.15cm

\noindent {\bf The case $(l_0,l_1,l_2,l_3)=(3,2,2,3)$.} 

\vskip.15cm

Now the arithmetic progression is given by $x_0^3,x_1^2,x_2^2,x_3^3$. A
possible approach would be to follow a similar argument as in the previous
case. That is, multiplying the formulas
$$
x_1^2=\frac{2x_0^3+x_3^3}{3},\ \ \ x_2^2=\frac{x_0^3+2x_3^3}{3}
$$
we get
$$
(3x_1x_2)^2=2x_0^6+5x_0^3x_3^3+2x_3^6.
$$
If $x_3=0$ then gcd$(x_2,x_3)=1$ yields $x_1^2=\pm 2$, a contradiction. So
we may suppose that $x_3\neq 0$, and we obtain
$$
Y^2=2X^6+5X^3+2
$$
with $X=x_0/x_3$ and $Y=3x_1x_2/x_3^3$. However, a calculation with MAGMA
gives that the rank of the Jacobian of the above hyperelliptic curve is two,
hence we cannot apply the classical Chabauty argument in this case. So we
follow a different method, which also makes it possible to exhibit an
elliptic curve (over some number field) having non-trivial Tate-Shafarevich
group.

For this purpose, observe that we have
$$
(-x_0x_3)^3=2d^2-(x_1x_2)^2,
$$ 
where $d$ denotes the increment of the progression. Using part i) of our Lemma we get
that there are two possible parametrizations given by
$$
x_1x_2=\pm(x^3+6x^2y+6xy^2+4y^3),\ d=\pm(x^3+3x^2y+6xy^2+2y^3),\
x_0x_3=-x^2+2y^2
$$
or
$$
x_1x_2=\pm(x^3+6xy^2),\ d=\pm(3x^2y+2y^3),\ x_0x_3=x^2-2y^2.
$$
Therefore from $x_1^2+d=x_2^2$ either
$$
x_1^4+dx_1^2-(x^3+6x^2y+6xy^2+4y^3)^2=0\tag{5}
$$
or
$$
x_1^4+dx_1^2-(x^3+6xy^2)^2=0\tag{6}
$$
follows, respectively. In the first case, the left hand side of
\thetag{5} can be considered as a polynomial of degree two in
$x_1^2$. Hence its discriminant must be a perfect square in ${\Bbb Z}$,
and we get the equation
$$
5x^6+54x^5y+213x^4y^2+360x^3y^3+384x^2y^4+216xy^5+68y^6=z^2
$$
in integers $x,y,z$. A simple calculation with MAGMA shows that the
Jacobian of the corresponding hyperelliptic curve
$$
Y^2=5X^6+54X^5+213X^4+360X^3+384X^2+216X+68
$$
is of rank zero (anyway it has three torsion points), and there is no rational point on the curve at all. Hence
in this case we are done. It is interesting to note, however, that this curve
does have points everywhere locally. We really do need this global information
on the rank of its Jacobian in order to decide it does not have any rational
points.

In case of \thetag{6} by a similar argument we obtain that
$d^2+4(x^3+6xy^2)^2=z^2$, whence
$$
4x^6+57x^4y^2+156x^2y^4+4y^6=z^2
$$
with certain integers $x,y,z$. Observe that $y=0$ yields a non-primitive solution.
Hence after putting $Y=z/2y^3$ and $X=x/y$, we get that to solve
the above equation it is sufficient to find all rational points on
the curve
$$
{\Cal C}_4:Y^2=f(X)=X^6+(57/4)X^4+39X^2+1.
$$
We show that the rational points on ${\Cal C}_4$ all have $X\in\{0,\infty\}$.

A straightforward computation shows that the rank of the Jacobian
${\Cal J}({\Bbb Q})$ of ${\Cal C}_4$ is two, so we cannot apply Chabauty's
method as before (cf. also \cite{CF}). We use part of the $2$-coverings of
${\Cal C}_4$ following \cite{B}. For details, see \cite{notes}. Let
$$
K=\Bbb{Q}(\alpha)=\Bbb{Q}[X]/(X^3+(57/4)X^2+39X+1).
$$
Over this field, we have
$$
f(X)=Q(X)R(X)=(X^2-\alpha)(X^4+(\alpha+57/4)X^2+\alpha^2+(57/4)\alpha+39).
$$
One easily gets that $\text{Res}(Q,R)$ is a unit outside
$S=\{\text{places}\ \frak{p}\ \text{of}\ K\ \text{dividing}\ 6\ \text{or}\
\infty\}$. Therefore, if $(X,Y)\in {\Cal C}_4({\Bbb Q})$ then we have
$$
\align
D_\delta&:~~ (Y_1)^2=\delta R(X)\\
L_\delta&:~~ (Y_2)^2=\delta Q(X)
\endalign
$$
for some $Y_1,Y_2\in K$ and $\delta\in K^*$ representing some element of the
finite group
$$
K(S,2):=\{[d]\in K^*/K^{*2} : 2\mid \text{ord}_\frak{p}(d)\ \text{for all
places}\ \frak{p}\notin S\}.
$$
Furthermore, since $N_{K[X]/{\Bbb Q}[X]}(Q)=f$, we see that
$N_{K/{\Bbb Q}}(\delta)\in {\Bbb Q}^{*2}$. Running through these finitely
many candidates, we see that the only class for which $D_\delta$ has points
locally at the places of $K$ above $2$ and $\infty$ is represented by
$\delta=1$. Over $K$, the curve $D_1$ is isomorphic to
$$
E:v^2=u^3-{{4\alpha+57}\over {2}}u^2-{{48\alpha^2+456\alpha-753}\over{16}}u,
$$
where $X=v/(2u)$. This curve has full $2$-torsion over $K$ and a full
$2$-descent or any $2$-isogeny descent gives a rank bound of two for $E(K)$.
However, one of the isogenous curves,
$$
E':Y^2=X^3+(4\alpha+57)X^2+(16\alpha^2+228\alpha+624)X
$$
has $S^{(2)}(E'/K)\simeq {\Bbb Z}/2{\Bbb Z}$, which shows that $E'(K)$ is
of rank zero, since $E'$ has $4$-torsion over $K$. This shows that $E$ has
non-trivial $2$-torsion in its Tate-Shafarevich group and that $E(K)$
consists entirely of torsion. In fact,
$$
E(K)=\{\infty,(0,0),((12\alpha^2+195\alpha+858)/32,0),
((-12\alpha^2-131\alpha+54)/32,0)\}.
$$
It follows that
$$
X({\Cal C}_4({\Bbb Q}))\subset X(D_1(K))=\{0,\infty\},
$$
where $X(.)$ denotes the set of the $X$-coordinates of the appropriate
points on the corresponding curve. This proves that for all the rational
points on ${\Cal C}_4$ we have $X\in\{0,\infty\}$, which implies the
theorem also in this case.

\vskip.15cm

\noindent
{\bf The cases $(l_0,l_1,l_2,l_3)=(2,2,3,3)$ and $(3,3,2,2)$.} 

\vskip.15cm

Again by symmetry, we may assume that $(l_0,l_1,l_2,l_3)=(2,2,3,3)$. Then
the progression is $x_0^2,x_1^2,x_2^3,x_3^3$, whence
$$
x_1^2=2x_2^3-x_3^3\ \ \ \text{and}\ \ \ x_0^2=3x_2^3-2x_3^3.
$$
If $x_3=0$ then the coprimality of $x_2$ and $x_3$ gives $x_1^2=\pm 2$,
which is a contradiction. Hence we may assume that $x_3\neq 0$, and we get
the equation
$$
y^2=F(x)=6x^6-7x^3+2
$$
with $x=x_2/x_3$, $y=x_0x_1/x_3^3$. Put $K={\Bbb Q}(\alpha)$ with
$\alpha=\root{3}\of{2}$ and observe that we have the factorization
$F(x)=G(x)H(x)$ over $K$ where
$$
G(x)=3\alpha x^4-3x^3-2\alpha x+2\ \ \ \text{and}\ \ \
H(x)=\alpha^2x^2+\alpha x+1.
$$
A simple calculation by MAGMA gives that Res$(G,H)$ is a unit
outside the set $S=\{\text{places}\ \frak{p}\ \text{of}\ K\ \text{dividing}\ 6\
\text{or}\ \infty\}$. Hence we can write
$$
3\alpha x^4-3x^3-2\alpha x+2=\delta z^2
$$
with some $z$ from $K$ and $\delta$ from the integers of $K$ dividing
$6$. Moreover, observe that the norm of $\delta$ is a square in
${\Bbb Z}$. Using that $\alpha-1$ is a fundamental unit of $K$, $2=\alpha^3$
and $3=(\alpha-1)(\alpha+1)^3$, local considerations show that
we can only have solutions with $x\in\QQ$ with both $G(x)$ and $H(x)\in K^{*2}$ 
if, up to squares, $\delta=\alpha-1$. We consider
$$
3\alpha x^4-3x^3-2\alpha x+2=(\alpha-1)z^2
$$
with $x\in{\Bbb Q}$ and $z\in K$. Now by the help of the point $(1,1)$, we can
transform this curve to Weierstrass form
$$
E: X^3+(-72\alpha^2-90\alpha-108)X+(504\alpha^2+630\alpha+798)=Y^2.
$$
We have $E(K)\simeq \ZZ$ as an Abelian group and the point
$(X,Y)=(-\alpha^2-1,12\alpha^2+15\alpha+19)$ is a non-trivial point on this
curve.
Again applying elliptic Chabauty with $p=5$, we get that the only
solutions of our original equation is $(x,z)=(1,1)$. Hence the theorem
follows also in this case.

\vskip.15cm

\noindent
{\bf The case $(l_0,l_1,l_2,l_3)=(2,3,3,2)$.}

\vskip.15cm

Now we have a progression $x_0^2,x_1^3,x_2^3,x_3^2$, and we can write
$$
x_0^2=2x_1^3-x_2^3\ \ \ \text{and}\ \ \ x_3^2=-x_1^3+2x_2^3.
$$
If $x_2=0$ then the coprimality of $x_1$ and $x_2$ gives $x_0^2=\pm 2$,
which is a contradiction. Hence we may assume that $x_2\neq 0$, and we are
led to the equation
$$
y^2=F(x)=-2x^6+5x^3-2
$$
with $x=x_1/x_2$, $y=x_0x_3/x_2^3$. Now we have the factorization
$F(x)=G(x)H(x)$ over $K={\Bbb Q}(\alpha)$ with $\alpha=\root{3}\of{2}$, where
$$
G(x)=\alpha^2x^4+(\alpha+2)x^3+(\alpha^2+2\alpha+1)x^2+(\alpha+2)x+\alpha^2
$$
and
$$
H(x)=-\alpha x^2+(\alpha^2+1)x-\alpha.
$$
One can easily verify that Res$(G,H)=1$. Thus we
obtain
$$
\alpha^2x^4+(\alpha+2)x^3+(\alpha^2+2\alpha+1)x^2+(\alpha+2)x+\alpha^2=
\delta z^2
$$
where $z\in K$ and $\delta$ is a unit of $K$. Moreover, as the norm of
$\delta$ is a square in ${\Bbb Z}$, we get that, up to squares, $\delta=1$ or $\alpha-1$. 
The case when $\delta=1$ yields the equation
$$
\alpha^2x^4+(\alpha+2)x^3+(\alpha^2+2\alpha+1)x^2+(\alpha+2)x+\alpha^2
=z^2
$$
in $x\in{\Bbb Q}$ and $z\in K$. We can transform this equation to an elliptic
one by the help of its point $(1,\alpha^2+\alpha+1)$. Then applying elliptic
Chabauty, the procedure ``Chabauty'' of MAGMA with $p=5$ in this
case gives that this equation has four solutions with $x\in{\Bbb Q}$, namely
$(x,z)=(0,1),(1,0),(\pm 1,1)$. Lifting these solutions to the original
problem, our theorem follows also in this case.

When $\delta=\alpha-1$, using $x=x_1/x_2$ we get the equation
$$
\alpha^2x_1^4+(\alpha+2)x_1^3x_2+(\alpha^2+2\alpha+1)x_1^2x_2^2
+(\alpha+2)x_1x_2^3+\alpha^2x_2^4=(\alpha-1)\gamma^2
$$
with some integer $\gamma$ of $K$. Writing now $\gamma$ in the form
$\gamma=u+\alpha v+\alpha^2 w$ with some $u,v,w\in{\Bbb Z}$ and comparing the
coefficients of $1$ and $\alpha$ in the above equation, a simple calculation
shows that $x_1^3x_2+x_1^2x_2^2+x_1x_2^3$ must be even. However, then
$2\mid x_1x_2$, and considering the progression $x_0^2,x_1^3,x_2^3,x_3^2$
modulo $4$ we get a contradiction. Hence the theorem follows also in this
case.

\vskip.15cm

\noindent
{\bf The case $(l_0,l_1,l_2,l_3)=(3,3,2,3)$ and $(3,2,3,3)$.}

\vskip.15cm

As previously, without loss of generality we may assume that
$(l_0,l_1,l_2,l_3)=(3,3,2,3)$. Then the progression is of the form
$x_0^3,x_1^3,x_2^2,x_3^3$. We note that using the cubes one would find
$3x_1^3=x_3^3+2x_0^3$ which leads to an elliptic curve. However, this
elliptic curve has positive rank, hence this approach does not work.

So we use some other argument. We have $x_1^3+x_3^3=2x_2^2$, whence
$$
x_1+x_3=2su^2,\ \ \ x_1^2-x_1x_3+x_3^2=sv^2,
$$
where $u,v,s\in{\Bbb Z}$ with $s\mid 3$. By considerations modulo $3$ we
obtain that only $s=1$ is possible. Hence $(2x_1-x_3)^2+3x_3^2=(2v)^2$ and
from part viii) of our Lemma we get that
$$
f(x,y):=3x^6+18x^5y+9x^4y^2-148x^3y^3-27x^2y^4+162xy^5-81y^6=2(\pm 4x_0)^3
\tag{7}
$$
in coprime integers $x,y$.

Note that the equation $f(x,y)=2z^3$ is invariant under the transformation
$(x,y,z)\mapsto(-3y,x,-3z)$. The two obvious solutions $(x,y,z)=(1,-1,-4)$ and
$(x,y,z)=(3,1,12)$ are interchanged by this involution.
 
We have the factorization $f(x,y)=g(x,y)h(x,y)$ with
$$
g(x,y)=(\alpha^2+2\alpha+1)x^3+(-2\alpha^3-\alpha^2+2\alpha+1)x^2y+
$$
$$
(3\alpha^2-26\alpha-13)xy^2+(-6\alpha^3-3\alpha^2+6\alpha+3)y^3
$$
and
$$
h(x,y)=(2\alpha^3+3\alpha^2-2\alpha+9)x^3
+(12\alpha^3+17\alpha^2-10\alpha+53)x^2y+
$$
$$
(6\alpha^3+9\alpha^2-6\alpha+27)xy^2+(-92\alpha^3-141\alpha^2+66\alpha-401)y^3
$$
over the
number field ${\Bbb Q}(\alpha)$ defined by a root $\alpha$ of the polynomial
$X^4+2X^3+4X+2$.

Using the same reasoning as before, we have that a rational solution to
$f(x,y)=2z^3$ with $x,y,z$ not all $0$, yields a solution to the system of
equations
$$\align
g(x,y)=&\delta (u_0+u_1\alpha+u_2\alpha^2+u_3\alpha^3)^3\\
  h(x,y)=&2/\delta (v_0+v_1\alpha+v_2\alpha^2+v_3\alpha^3)^3
  \endalign
$$
with $x,y,u_0,\ldots,v_3\in\QQ$ and where $\delta$ is a representative of an
element of the finite group $K(S,3)$, with
$S=\{\text{places}\ \frak{p}\ \text{of}\ K\ \text{dividing}\ 6\ \text{or}\
\infty\}$.
For each $\delta$, the equations above can be expressed as eight homogeneous
equations of degree $3$, describing some non-singular curve in $\PP^8$ over
$\QQ$. The only values of $\delta$ for which this curve is locally solvable at
$3$ are
$$\delta_1=(\alpha^3+2\alpha^2-2\alpha-2)/2\text{ and }
\delta_2=(\alpha^3+4\alpha^2+6\alpha+2)/2.$$
These values correspond to the obvious solutions with $(x,y)=(1,-1)$ and
$(x,y)=(3,1)$ respectively.

We now determine the $K$-rational points on the curve
$$g(x,y)=\delta_1 z_1^3$$
with $x/y\in\QQ$. Using the $K$-rational point $(x:y:z)=(1:-1:-2\alpha)$, we can see that
this curve is isomorphic to the elliptic curve
$$E: Y^2=X^3-48\alpha^3+33\alpha^2+480\alpha+210.$$
Using a $2$-descent we can verify that $E(K)$ has rank at most $3$ and some
further computations show that $E(K)\simeq\ZZ^3$, where the points with
$X$-coordinates
$$\align
    (-2\alpha^3 + 13\alpha^2 - 28\alpha + 44)/9,&\\
    (16\alpha^3 + 52\alpha^2 + 14\alpha - 1)/9,&\\
    (2\alpha^3 + 3\alpha^2 - 14\alpha - 6)/3&
\endalign$$
generate a finite index subgroup with index prime to $6$.
The function $x/y$ on the curve $g(x,y)=\delta_1 z_1^3$ yields a degree $3$
function on $E$ as well.

Using the Chabauty-method described in \cite{B} and implemented in MAGMA 2.11
as {\tt Chabauty}, using $p=101$, we determine that the given point is in fact
the only one with $x/y\in\QQ$. For details, see \cite{notes}.

For $\delta_2$ we simply observe that using the involution
$(x,y)\mapsto(-3y,x)$, we can reduce this case to the computations we have
already done for $\delta_1$.

We conclude that $(x,y)=(1,-1)$ and $(x,y)=(3,1)$ give the only solutions to
$f(x,y)=2z^3$. These solutions correspond to the arithmetic progressions
$(0,1,2,3)$ (which up to powers of $2,3$ indeed consists of second and third
powers), $(1,1,1,1)$ and their $\ZZ_{\{2,3\}}^*$-equivalent counterparts.
\qed
\enddemo

\heading{Acknowledgements}
\endheading

The authors are grateful to the referee for his useful and helpful remarks.

\enddocument